\documentclass{ifacconf}
\usepackage{graphicx}      
\usepackage{natbib}        
\usepackage[utf8]{inputenc}

\usepackage{amsmath}

\usepackage{xspace}

\def\fins{p}

\begin{document}
\begin{frontmatter}

\title{Discrete-time Flatness-based Controller Design using an Implicit Euler-discretization} 

\thanks[footnoteinfo]{ This work has been supported by the Austrian Science Fund (FWF) under grant number P 32151.}

\author[First]{Johannes Diwold} 
\author[Second]{Bernd Kolar} 
\author[First]{Markus Schöberl}

\address[First]{Institute of Automatic Control and Control Systems Technology,	Johannes Kepler University Linz, Altenbergerstraße 66, 4040 Linz, Austria (e-mail: johannes.diwold@jku.at, markus.schoeberl@jku.at)}
\address[Second]{Magna Powertrain Engineering Center Steyr GmbH \& Co KG, Steyrer	Str. 32, 4300 St. Valentin, Austria (e-mail: bernd\_kolar@ifac-mail.org)}

\begin{abstract}
In this contribution, we present a constructive method to derive flat sampled-data models for continuous-time flat systems through an implicit Euler-discretization. We show how the sampled-data model can be used subsequently for a flatness-based controller design, and illustrate our results with the well-known planar VTOL aircraft example.
\end{abstract}


\end{frontmatter}

\section{Introduction}

In the 1990s, the concept of flatness has been introduced by Fliess,
L\'{e}vine, Martin and Rouchon for nonlinear continuous-time systems
(see e.g. \cite{FliessLevineMartinRouchon:1995} and \cite{FliessLevineMartinRouchon:1999}).
A nonlinear continuous-time system
\begin{equation}
\dot{x}=f(x,u)\label{eq:system_equation_continuous}
\end{equation}
with $\dim(x)=n$ states, $\dim(u)=m$ inputs and smooth functions
$f$ is flat, if there exists a one-to-one correspondence between
solutions $(x(t),u(t))$ of \eqref{eq:system_equation_continuous}
and solutions $y(t)$ of a trivial system (sufficiently smooth but
otherwise arbitrary trajectories). This one-to-one correspondence
implies that for flat systems there exist flat outputs 
\begin{equation}
y=\varphi(x,u,\dot{u},\dots,u^{(q)})
\end{equation}
such that all states $x$ and inputs $u$ can be parameterized by
$y$ and its time derivatives, i.e., 
\begin{equation}
(x,u)=F(y,\dot{y},\dots,y^{(r)})\,.\label{eq:continuous_time_parameterization}
\end{equation}
We call $F$ the parameterizing map with respect to the flat output
$y$. The popularity of flatness is due to the fact that it allows
an elegant solution to motion planning problems as well as the systematic
design of tracking controllers. For the practical implementation of
such a tracking controller,
it is
necessary to evaluate the continuous-time control law at
a sufficiently high sampling rate. If, however, the processing unit or actuators/sensors are limited to lower sampling rates, a discrete
evaluation of a continuous-time control law may lead to unsatisfactory
results. As known from linear systems theory, an appropriate alternative
is to design the controller directly for a suitable discretization.
Motivated by these considerations, in \cite{DiwoldKolarSchoeberl2022-1}
we employed the concept of discrete-time flatness (\cite{KaldmaeKotta:2013},
\cite{GuillotMillerioux2020} and \cite{DiwoldKolarSchoeberl2022})
and designed a discrete-time flatness-based tracking control for an
explicit Euler-discretization of a gantry crane. Experiments with a
laboratory setup have shown that the discrete-time approach is indeed
more robust with respect to large sampling times than its continuous-time
counterpart. However, also with the discrete-time approach, the performance
decreases for large sampling times. Therefore, the objective of the
present contribution is to replace the explicit Euler-discretization
used in \cite{DiwoldKolarSchoeberl2022-1} by an implicit Euler-discretization,
which is known to be preferable at least in terms of numerical aspects.

In general, the flatness property is not preserved by an exact or
approximate discretization. Thus, in \cite{DiwoldKolarSchoeberl2022-1}
we have proposed a constructive method that combines a suitable state
transformation into a structurally flat triangular form and a subsequent explicit Euler-discretization such that the flatness of the system is preserved. Although only
constructive, the method is applicable for a large class of nonlinear
systems, e.g. the induction motor, the VTOL aircraft, and the unmanned
aerial vehicle discussed in \cite{Greeff2021}. In this paper, we
adapt the method of \cite{DiwoldKolarSchoeberl2022-1} by using an implicit Euler-discretization.
Accordingly, we transform the system \eqref{eq:system_equation_continuous} into a structurally flat triangular form first, and perform an implicit Euler-discretization afterwards. Since the implicit Euler-discretization preserves the triangular structure, it also preserves the flatness of the system. Moreover, the triangular structure of the discrete-time system allows to compute the parameterizing map in a systematic way. Once the parameterization has been determined, the design of a tracking control law based on an exact linearization by a dynamic feedback is a straightforward task, see \cite{DiwoldKolarSchoeberl2022}.

The paper is organized as follows: In Section \ref{sec:Flatness-of-Discrete-time}
we recapitulate the notion of flatness for discrete-time systems
as proposed in \cite{DiwoldKolarSchoeberl2022}. Next, in Section \ref{sec:ImpEulerContDes}
we show that the implicit Euler-discretization of a continuous-time system in a structurally flat triangular form preserves the flatness of the system. Furthermore, we show how the parameterization can be determined in a systematic way, and how it can be used for the controller design. In Section
\ref{sec:example}, we illustrate the proposed method
by the practical example of a VTOL aircraft, design a discrete-time flatness-based tracking controller, and present simulation results.
%
\section{Flatness of Discrete-time Systems}
\label{sec:Flatness-of-Discrete-time}
In this section we briefly recall the notion of discrete-time flatness as proposed in \cite{DiwoldKolarSchoeberl2022}. 

We consider discrete-time systems in state representation
\begin{equation}
\label{eq:system_equations_disc}
x^{i,+}=f^i(x,u)\,,\phantom{a}i=1,\dots,n
\end{equation}
with $\dim(x)=n$, $\dim(u)=m$ and smooth functions $f^{i}(x,u)$.
Furthermore, we assume that the system meets the submersivity condition, i.e. that the Jacobian matrix of $f$ with respect to $(x,u)$ meets
\begin{equation}
\textrm{rank}(\partial_{(x,u)}f)=n\,.
\label{eq:submersivity_condition}
\end{equation}
This condition is necessary for reachability and consequently also for flatness. As stated in \cite{DiwoldKolarSchoeberl2022}, flatness for discrete-time systems can be characterized by a one-to-one correspondence between the system trajectories $(x(k),u(k))$ and trajectories $y(k)$ of a trivial system, i.e., arbitrary trajectories that need not satisfy any difference equation. This one-to-one correspondence can be expressed in terms of two maps
\begin{equation}
\label{eq:parm_map_time_variant}
(x(k),u(k))=F(y(k-r_1),\dots,y(k+r_2))
\end{equation}
and
\begin{equation}
\label{eq:flat_output_time_variant}
\begin{aligned}
y(k)=\varphi(x(k-q_1),u(k-q_1)\dots,x(k+q_2),u(k+q_2))
\end{aligned}
\end{equation}
with some integers $r_1,r_2,q_1,q_2\geq0$. As explained in \cite{DiwoldKolarSchoeberl2022}, the composition of \eqref{eq:parm_map_time_variant} with the occurring shifts of \eqref{eq:flat_output_time_variant} must yield the identity map, and vice versa. Furthermore, since the trajectories $y(k)$ of a trivial system are arbitrary, substituting \eqref{eq:parm_map_time_variant} into \eqref{eq:system_equations_disc} must also yield an identity. The map \eqref{eq:flat_output_time_variant} can be further simplified by taking into account that every trajectory $\dots,x(k-1),u(k-1),x(k),u(k),x(k+1),u(k+1),\dots$ satisfies \eqref{eq:system_equations_disc}. Hence, the forward-shifts of $x(k)$ can be expressed as functions of $x(k),u(k),u(k+1),\dots$ via the successive compositions
\begin{equation}
\label{eq:to_replace_forward_shifts_of_x}
\begin{aligned}
x(k+1)&=f(x(k),u(k))\\
x(k+2)&=f(f(x(k),u(k)),u(k+1))\,.\\
&\phantom{a}\vdots
\end{aligned}
\end{equation}
A similar argument holds for the backward-direction. Since \eqref{eq:system_equations_disc} meets the submersivity condition \eqref{eq:submersivity_condition}, there always exist $m$ functions $g(x,u)$ such that the map
\begin{equation}
\label{eq:extending_f_to_diffeomorphism}
x^+=f(x,u)\,,\phantom{aa}\zeta=g(x,u)
\end{equation}
is locally invertible.
If we denote by $(x,u)=\psi(x^+,\zeta)$ its inverse
\begin{align}
\label{eq:inverted_extension}
\begin{aligned}x & =\psi(x^+,\zeta)\,,\end{aligned}
& u=\psi(x^+,\zeta)\,,
\end{align}
then we can express the backward-shifts of $(x(k),u(k))$ as functions of $x(k),\zeta(k-1),\zeta(k-2),\dots$ via the successive compositions
\begin{equation}
\label{eq:to_replace_backward_shifts_of_x_and_u}
\begin{aligned}
(x(k-1),u(k-1))&=\psi(x(k),\zeta(k-1))\\
(x(k-2),u(k-2))&=\psi(\psi_{x}(x(k-1),\zeta(k-1)),\zeta(k-2))\,.\\
&\phantom{a}\vdots
\end{aligned}
\end{equation}
Thus, with \eqref{eq:to_replace_forward_shifts_of_x} and \eqref{eq:to_replace_backward_shifts_of_x_and_u}, the map \eqref{eq:flat_output_time_variant} can be written as
\begin{equation}
\label{eq:flat_output_as_functions_of_zeta_and_u}
y(k)=\varphi(\zeta(k-q_1),\dots,\zeta(k-1),x(k),u(k),\dots,u(k+q_2))\,.
\end{equation}
The map \eqref{eq:flat_output_as_functions_of_zeta_and_u} is called a flat output of system the \eqref{eq:system_equations_disc}, and the corresponding map \eqref{eq:parm_map_time_variant} describes the parameterization of the system variables $(x(k),u(k))$ by this flat output.
\begin{rem}
	\label{rem:choice_zeta}The flatness
	of the system (\ref{eq:system_equations_disc}) does not depend on the choice
	of the functions $g(x,u)$ of \eqref{eq:extending_f_to_diffeomorphism}. The representation \eqref{eq:flat_output_as_functions_of_zeta_and_u}
	of the flat output may differ, but the parameterization \eqref{eq:parm_map_time_variant}
	is not affected. 
\end{rem}
Note also that considering backward-shifts in both \eqref{eq:parm_map_time_variant} and \eqref{eq:flat_output_as_functions_of_zeta_and_u} is actually not necessary. If there exist a parameterizing map \eqref{eq:parm_map_time_variant} and a flat output \eqref{eq:flat_output_as_functions_of_zeta_and_u}, then one can always define a new flat output as the $r_1$-th backward-shift of the original flat output. The corresponding parameterizing map is then of the form
\begin{equation}
\label{eq:parameterizing_map_forward_shifts_only}
(x(k),u(k))=F(y(k),\dots,y(k+r))
\end{equation}
with $r=r_1+r_2$ and does not contain backward-shifts. Similarly, we may define a new flat output as the $q_{1}$-th forward-shift
of the original flat output. The new flat output reads as
\begin{equation}
\label{eq:ch3_flat_output_trajectory}
y(k)=\varphi(x(k),u(k),\dots,u(k+q))
\end{equation}
with $q=q_{1}+q_{2}$,
and the corresponding parameterizing map is of the form \eqref{eq:parm_map_time_variant}.
\begin{rem}
Systems that possess both a flat output \eqref{eq:ch3_flat_output_trajectory} independent of backward-shifts as well as a parameterizing map \eqref{eq:parameterizing_map_forward_shifts_only} independent of backward-shifts represent a special subset of flat systems which will be denoted in the following as forward-flat. This class of systems is studied in detail e.g. in \cite{KaldmaeKotta:2013}, \cite{KolarSchoberlDiwold:2019}, and \cite{KolarDiwoldSchoberl:2019}.
\end{rem}
In the sequel, for the parameterizing map \eqref{eq:parm_map_time_variant} we write
\begin{equation}
(x,u)=F(y_{[-R_1]},\dots,y_{[R_2]})\,,
\label{eq:parameterizing_map_forward_and_backward_shifts}
\end{equation}
i.e., we remove $k$ and use subscripts in brackets\footnote{For a single forward-shift we also use the superscript $+$.} to denote the forward- and backward-shifts instead. The multi-index $R_1=(r_1^1,\dots,r_1^m)$ denotes the highest backward-shift of the flat output $(y^1,\dots,y^m)$ that occur in the parameterizing map. Likewise the multi-index $R_2$ denotes the highest forward-shifts. For the flat output we write
\begin{equation}
\label{eq:flat_output_disc_zta}
y=\varphi(\zeta_{[-q_1]},\dots,\zeta_{[-1]},x,u,u_{[1]},\dots,u_{[q_2]})\,.
\end{equation}
For a concise geometric definition of discrete-time flatness we refer to \cite{DiwoldKolarSchoeberl2022}.
\section{Implicit Euler-Discretization and Controller Design}

\label{sec:ImpEulerContDes}For the design of discrete-time flatness-based
controllers, a sampled-data model of the plant has to be derived.
Preferably an exact discretization is used, however, this is a difficult
task for most nonlinear systems. Furthermore, an exact discretization
does not preserve the flatness of a system in general, see \cite{DiwoldKolarSchoeberl2022}.
From a practical point of view, using an Euler-discretization is much
more convenient. However, even the Euler-discretization does not necessarily
preserve flatness. In \cite{DiwoldKolarSchoeberl2022-1}, we have
shown that applying a suitable state transformation before performing
an explicit Euler-discretization is useful for deriving a forward-flat
sampled-data system. In this section, we adapt the method of \cite{DiwoldKolarSchoeberl2022-1} to use an implicit Euler-discretization instead, which is known to be preferable at least in terms of numerical aspects.
As we will see below, in contrast to the procedure proposed in \cite{DiwoldKolarSchoeberl2022-1},
the resulting discretization is flat but not necessarily forward-flat.
\subsection{Explicit and Implicit Euler-Discretization\label{subsec:Implicit-Euler-Discretization}}

In general, the exact discretization of a time-invariant system \eqref{eq:system_equation_continuous}
can be derived by evaluating
\[
x^{+}=x+\int_{0}^{T_{s}}f(x(t),u)\mathrm{d}t\,,
\]
where the input $u$ is constant during the sampling interval. For
most nonlinear systems, evaluating the integral $\int_{0}^{T_{s}}f(x(t),u)\mathrm{d}t$
is a difficult task, and hence a suitable approximation is used instead.
The explicit Euler-discretization $x^{+}=x+T_{s}f(x,u)$ results from
the approximation
\[
\int_{0}^{T_{s}}f(x(t),u)\mathrm{d}t\approx T_{s}f(x,u),
\]
i.e., the integral is approximated by the area of a rectangle determined
by the left-sided value of $x$. Contrary, the implicit Euler-discretization
follows from the approximation
\[
\int_{0}^{T_{s}}f(x(t),u)\mathrm{d}t\approx T_{s}f(x^{+},u),
\]
where the right-sided value of $x$ is used. In contrast to the explicit
Euler-discretization, the implicit Euler-discretization
\begin{equation}
x^{+}=x+T_{s}f(x^{+},u)\label{eq:implicit_euler_discretization_implicit_form}
\end{equation}
does not immediately result in a state representation \eqref{eq:system_equations_disc}.
Thus, in order to apply the notion of discrete-time flatness of Section \ref{sec:Flatness-of-Discrete-time},
we need to show that there exists a state representation for \eqref{eq:implicit_euler_discretization_implicit_form},
at least locally.
\begin{thm}
	\label{thm:For-a-sufficiently}Consider the implicit Euler-discretization
	\eqref{eq:implicit_euler_discretization_implicit_form} of a nonlinear
	continuous-time system \eqref{eq:system_equation_continuous} with
	an equilibrium $(x_{s},u_{s})$. For a sufficiently small sampling
	time $T_{s}>0$, there exists a state representation
	\begin{equation}
	x^{+}=\tilde{f}(x,u)\label{eq:implicit_euler_discretization_explicit_form}
	\end{equation}
	in a neighborhood of $(x_{s},u_{s})$.
\end{thm}
The proof is based on an application of the implicit function theorem (see e.g. \cite{NijmeijervanderSchaft:1990}) around the point $(x^+,x,u,T_s)=(x_s,x_s,u_s,0)$.
Although Theorem \ref{thm:For-a-sufficiently} guarantees the existence
of a state representation \eqref{eq:implicit_euler_discretization_explicit_form} for small sampling times $T_s>0$,
a symbolic expression for $\tilde{f}(x,u)$ may be hard to find. However, in the remainder of this paper a symbolic expression for $\tilde{f}(x,u)$ will not be needed. 
\subsection{Implicit Euler-discretization of Triangular Forms}
\label{subsec:triangular}
As in \cite{DiwoldKolarSchoeberl2022-1}, let us now consider a flat continuous-time system which is decomposed into $p$ subsystems 
\begin{align}
\begin{aligned}\dot{{x}}_{\fins} & =f_{\fins}({x}_{\fins},{x}_{\fins-1},{u}_{\fins-1})\\
\dot{{x}}_{\fins-1} & =f_{\fins-1}({x}_{\fins},{x}_{\fins-1},{x}_{\fins-2},{u}_{\fins-1},{u}_{\fins-2})\\
& \phantom{a}\vdots\\
\dot{{x}}_{2} & =f_{2}({x}_{\fins},\dots,{x}_{1},{u}_{\fins-1},\dots,{u}_{1})\\
\dot{{x}}_{1} & =f_{1}({x}_{\fins},\dots,{x}_{1},{u}_{\fins-1},\dots,{u}_{1},{u}_{0})
\end{aligned}
\label{eq:structurally_flat_triangular_form_continuous_time_system}
\end{align}
with the state
\begin{equation}
\label{eq:def_x_blocks_struct_flat_triangular_form_cont}
\begin{aligned}
{x} & =({x}_{\fins},{x}_{\fins-1}\dots,{x}_{1})\\
& =(y_{\fins},(y_{\fins-1},\hat{x}_{\fins-1}),\dots,(y_{1},\hat{x}_{1}))
\end{aligned}
\end{equation}
and the input ${u}=({u}_{\fins-1},\dots,{u}_{0})$ (the components $y_{k}$ and/or ${u}_{k}$ might also be empty). We furthermore assume that the system \eqref{eq:structurally_flat_triangular_form_continuous_time_system} meets the rank conditions
\begin{align}
\begin{aligned}\mathrm{rank}(\partial_{{u}_{0}}f_{1}) & =\dim(f_{1})=\dim({u}_{0})\\
\textrm{\ensuremath{\mathrm{rank}}}(\partial_{(\hat{x}_{k},{u}_{k})}f_{k+1}) & =\dim(f_{k+1})=\dim((\hat{x}_{k},{u}_{k}))
\end{aligned}
\label{eq:rank_cond}
\end{align}
for $k=1,\dots,\fins-1$. Note that the subscript $k$ denotes the corresponding subsystem $\dot{x}_k=f_k$ with $\dim(x_k)=\dim(f_k)\geq 1$. Accordingly, with the superscript of $f^i_k$ we denote the $i$-th component of $f_k$, likewise for the components of $x_k$ and $u_k$. The system representation \eqref{eq:structurally_flat_triangular_form_continuous_time_system}, which we also refer to as a structurally flat triangular form, has the beneficial property that the one-to-one correspondence between solutions of \eqref{eq:structurally_flat_triangular_form_continuous_time_system} and solutions of a trivial system\footnote{Sufficiently smooth but otherwise arbitrary trajectories that need not satisfy any differential equation.} can be determined in a systematic way. For this purpose, we consider $y=(y_p,\dots,y_1)$ as smooth trajectories and compute the corresponding trajectories $(x,u)$ by solving the set of equations from top to bottom. Accordingly, we write the first subsystem as \begin{equation}
\label{eq:rewritten_first_subsystem}
\dot{y}_p=f_p(y_p,y_{p-1},\hat{x}_{p-1},u_{p-1})\,,
\end{equation}cf. \eqref{eq:def_x_blocks_struct_flat_triangular_form_cont}. Since the rank condition \eqref{eq:rank_cond} holds for $k=p-1$, the equations \eqref{eq:rewritten_first_subsystem} can be solved for $(\hat{x}_{p-1},u_{p-1})$ as functions of $y_p,\dot{y}_p,y_{p-1}$. Furthermore, the remaining states of $x_{p-1}$ correspond to $y_{p-1}$, and thus the trajectories $(x_{p-1},u_{p-1})$ of the first subsystem are completely parameterized by the trajectories $y_p,\dot{y}_p,y_{p-1}$. In a next step, due to the condition \eqref{eq:rank_cond} for $k=p-2$, we can likewise solve the set of equations of the second subsystem for $(\hat{x}_{p-2},u_{p-2})$, and parameterize the trajectories $(x_{p-2},u_{p-2})$ by the trajectories $y_p,\dot{y}_p,\ddot{y}_p,y_{p-1},\dot{y}_{p-1},y_{p-2}$. Continuing this procedure leads finally to a parameterizing map of the form \eqref{eq:continuous_time_parameterization} together with the flat output $y=(y_p,\dots,y_1)$ defined according to \eqref{eq:def_x_blocks_struct_flat_triangular_form_cont}.


Let us now apply the implicit Euler-discretization \eqref{eq:implicit_euler_discretization_implicit_form} to a system in the structurally flat triangular form \eqref{eq:structurally_flat_triangular_form_continuous_time_system}, i.e.
\begin{align}
\begin{aligned}{{x}}^+_{\fins} & =x_{p}+T_sf_{\fins}({x}^+_{\fins},{x}^+_{\fins-1},{u}_{\fins-1})\\
{{x}}^+_{\fins-1} & =x_{\fins-1}+T_sf_{\fins-1}({x}^+_{\fins},{x}^+_{\fins-1},{x}^+_{\fins-2},{u}_{\fins-1},{u}_{\fins-2})\\
& \phantom{a}\vdots\\
{{x}}^+_{2} & =x_{2}+T_sf_{2}({x}^+_{\fins},\dots,{x}^+_{1},{u}_{\fins-1},\dots,{u}_{1})\\
{{x}}^+_{1} & =x_1+T_sf_{1}({x}^+_{\fins},\dots,{x}^+_{1},{u}_{\fins-1},\dots,{u}_{1},{u}_{0})
\end{aligned}
\label{eq:structurally_flat_triangular_form_discrete_time_system}
\end{align}
together with \eqref{eq:def_x_blocks_struct_flat_triangular_form_cont} and
\begin{align}
\begin{aligned}\mathrm{rank}(\partial_{{u}_{0}}f_{1}) & =\dim(f_{1})=\dim({u}_{0})\\
\textrm{\ensuremath{\mathrm{rank}}}(\partial_{(\hat{x}_{k}^+,{u}_{k})}f_{k+1}) & =\dim(f_{k+1})=\dim((\hat{x}_{k}^+,{u}_{k}))\,.
\end{aligned}
\label{eq:rank_cond_disc}
\end{align}
It can be observed that the implicit Euler-discretization \eqref{eq:implicit_euler_discretization_implicit_form} preserves the triangular structure, and thus we can systematically determine the one-to-one correspondence between solutions, as in the continuous-time case. For this purpose, we consider $y=(y_p,\dots,y_1)$ as trajectories (sequence), and determine the corresponding trajectories $(x,u)$ as functions of forward- and backward-shifted trajectories $y$, also denoted by $\dots,y_{[-1]},y,y_{[1]},\dots$. In a first step, we write the topmost equations as 
\begin{equation}
y_{p}^+=y_{p}+T_sf_p(y_{p}^+,y_{p-1}^+,\hat{x}_{p-1}^+,u_{p-1})\,.
\label{eq:topmost_disc_sys}
\end{equation} 
Due to \eqref{eq:rank_cond_disc} for $k=p-1$, we can solve the equations \eqref{eq:topmost_disc_sys} for $(\hat{x}^+_{p-1},u_{p-1})$ as functions of $y_p,y_p^+,y_{p-1}$ (or accordingly $y_p,y_{p,[1]},y_{p-1}$). For the parameterization of the non-shifted state trajectories $\hat{x}_{p-1}$, however, we need an additional backward-shift, i.e. $\hat{x}_{p-1}$ is a function of $y_{p,[-1]},y_{p},y_{p-1,[-1]}$. Since the remaining states of $x_{p-1}$ correspond to the components $y_{p-1}$, the trajectories $(x_{p-1},u_{p-1})$ are completely parameterized by the trajectories $y_{p,[-1]},y_p,y_{p,[1]},y_{p-1,[-1]},y_{p-1}$. Continuing this procedure for the remaining equations leads finally to a parameterizing map of the form \eqref{eq:parameterizing_map_forward_and_backward_shifts} together with the flat output $y=(y_p,\dots,y_1)$ defined according to \eqref{eq:def_x_blocks_struct_flat_triangular_form_cont}.

Let us briefly summarize the main result of this section: Since the implicit Euler-discretization of \eqref{eq:structurally_flat_triangular_form_continuous_time_system} preserves the triangular structure, there exists also a one-to-one correspondence of solutions between the discrete-time system \eqref{eq:structurally_flat_triangular_form_discrete_time_system} and solutions of a trivial system $y$ (arbitrary trajectories). However, since flatness as defined in Section \ref{sec:Flatness-of-Discrete-time} is restricted to systems in state representation \eqref{eq:system_equations_disc}, we formulate the following lemma in terms of the explicit state representation of Theorem \ref{thm:For-a-sufficiently}.
\begin{lem}
	\label{lem:preserving_flatness_of_str_form}
 The implicit Euler-discretization \eqref{eq:implicit_euler_discretization_explicit_form} of a continuous-time system in a structurally flat triangular form \eqref{eq:structurally_flat_triangular_form_continuous_time_system} is flat. 
\end{lem} 
Subsequently, in order to be able to apply Lemma \ref{lem:preserving_flatness_of_str_form} to a larger class of systems, we also consider state- and input transformations resulting in a structurally flat triangular form \eqref{eq:structurally_flat_triangular_form_continuous_time_system}.
\subsection{State- and Input Transformations}
First of all we want to emphasize that in general state- and input transformations
\begin{subequations}
	\label{eq:state_input_transformation_sub}
	\begin{align}
	\bar{x}&=\Phi_x(x)
	\label{eq:state_transformation}\\
	\bar{u}&=\Phi_u(x,u)
	\label{eq:input_transformation}
	\end{align}
\end{subequations}
and the implicit Euler-discretization \eqref{eq:implicit_euler_discretization_implicit_form} do not commute. In other words, if two continuous-time systems \eqref{eq:system_equation_continuous} are related by a transformation \eqref{eq:state_input_transformation_sub}, then their implicit Euler-discretizations are not necessarily related by the same state- and input transformation. However, since the choice of a state is not unique, it is of course justified to perform a state transformation \eqref{eq:state_transformation} before a discretization. With input transformations \eqref{eq:input_transformation}, on the other hand, the situation is different. The input $u$ of the continuous-time system \eqref{eq:system_equation_continuous} might be physically motivated or determined by the available actuators. Thus, a sampled-data system with the original input $u$ is preferable for the controller design.
\begin{lem}
	\label{lem:input_one_to_one}
If two continuous-time
systems \eqref{eq:system_equation_continuous} are related by an input transformation \eqref{eq:input_transformation}, then there exists a one-to-one correspondence between the solutions of their implicit Euler-discretizations \eqref{eq:implicit_euler_discretization_implicit_form}.
\end{lem}
\begin{pf}
	The implicit Euler-discretization for the system \eqref{eq:system_equation_continuous}
	reads as
	\begin{equation}
	x^{+}=x+T_{s}f(x^{+},u)\,,\label{eq:implicit_euler_disc_without_input_transf}
	\end{equation}
	while for the system $\dot{x}=f(x,{\Phi}^{-1}_{u}(x,\bar{u}))$ it reads
	as
	\begin{equation}
	\begin{aligned}x^{+}=x+T_{s}f(x^{+},{\Phi}^{-1}_{u}(x^{+},\bar{u}))\,.\end{aligned}
	\label{eq:implicit_euler_disc_with_input_transf}
	\end{equation}
	It can be observed that the equations \eqref{eq:implicit_euler_disc_without_input_transf}
	and \eqref{eq:implicit_euler_disc_with_input_transf} are related
	via the map $u={\Phi}^{-1}_{u}(x^{+},\bar{u})$. Thus, while the state trajectories for both discretizations coincide, the corresponding input trajectories are related by the maps\footnote{Note that the maps $u={\Phi}^{-1}_{u}(x^{+},\bar{u})$ and $\bar{u}={\Phi}_{u}(x^{+},{u})$ are no input transformations, they only describe the relation between the trajectories of the state and the original as well as the tranformed input.} $u={\Phi}^{-1}_{u}(x^{+},\bar{u})$ and $\bar{u}={\Phi}_{u}(x^{+},{u})$.
\end{pf}
With Lemma \ref{lem:input_one_to_one} we can now state and prove the main result of this section.
\begin{thm}
	\label{thm:Discretization_of_structurally_flat_triangular_forms}If
	a continuous-time system \eqref{eq:system_equation_continuous} can
	be transformed into a structurally flat triangular form \eqref{eq:structurally_flat_triangular_form_continuous_time_system}
	via a state- and input transformation \eqref{eq:state_input_transformation_sub},
	then the implicit Euler-discretization of the system $\dot{\bar{x}}=f(\bar{x},u)$,
	obtained from \eqref{eq:system_equation_continuous} by applying only
	the state transformation \eqref{eq:state_transformation}, is flat.
\end{thm}
\begin{pf}
By assumption, after performing the state- and input transformation \eqref{eq:state_input_transformation_sub}, the system $\dot{\bar{{x}}}=f(\bar{x},\bar{u})$ is in a structurally flat triangular form. Due to Lemma \ref{lem:preserving_flatness_of_str_form} the implicit Euler-discretization of $\dot{\bar{{x}}}=f(\bar{x},\bar{u})$ is flat. Because of Lemma \ref{lem:input_one_to_one}, there exists a one-to-one correspondence between the solutions of the implicit Euler-discretizations of $\dot{\bar{{x}}}=f(\bar{x},\bar{u})$ and $\dot{\bar{{x}}}=f(\bar{x},{u})$. Thus, also the implicit Euler-discretization of $\dot{\bar{{x}}}=f(\bar{x},{u})$, obtained by applying only the state transformation \eqref{eq:state_transformation} to the system \eqref{eq:system_equation_continuous}, is flat.
\end{pf}
Hence, it is possible to derive a flat sampled-data representation with the original input $u$ of the continuous-time system \eqref{eq:system_equation_continuous}. For deriving the corresponding parameterizating map $(\bar{x},u)=F(y_{[-R_1]},\dots,y_{[R_2]})$ in an efficient way, however, it can be useful to compute additionally also the discretization of the system \eqref{eq:system_equation_continuous} after applying both the state- and input transformation \eqref{eq:state_input_transformation_sub}. Because of the triangular structure of the resulting system, the parameterization $(\bar{x},\bar{u})=F(y_{[-R_1]},\dots,y_{[R_2]})$ can be computed easily as it is explained in Section 3.2. The parameterization of the original input $u$ follows then immediately via the relation $u={\Phi}^{-1}_{u}(\bar{x}^{+},\bar{u})$. The obtained parameterizing map $(\bar{x},u)=F(y_{[-R_1]},\dots,y_{[R_2]})$ can then be used e.g. for the controller design which we discuss in the next section.

\begin{rem}
	In fact, many continuous-time flat systems possess a structurally
	flat triangular form \eqref{eq:structurally_flat_triangular_form_continuous_time_system}, e.g. the gantry crane discussed in \cite{DiwoldKolarSchoeberl2022-1}
	and the model of an unmanned aerial vehicle used in \cite{Greeff2021}.
	In \cite{GstoettnerKolarSchoeberl2020,GstoettnerKolarSchoeberl2021} or \cite{NicolauRespondek:2020},
	it is even shown how the state- and input transformations \eqref{eq:state_input_transformation_sub}
	resulting in a structurally flat triangular form can be determined
	systematically for a subclass of flat systems with two inputs.
\end{rem}

\subsection{Controller Design via Dynamic Feedback Linearization}
\label{subsec:Controller-Design}
In the following, we assume that the continuous-time flat system \eqref{eq:system_equation_continuous}
has been discretized by an implicit Euler-discretization following
the procedure discussed above.
Although the corresponding explicit system equations
\begin{equation}
\bar{x}^{+}=\tilde{f}(\bar{x},u)\label{eq:state_representation_controller_design}
\end{equation}
may be symbolically unknown, a flatness-based control law can be derived in two steps. First, we derive a dynamic feedback that would exactly linearize the system \eqref{eq:state_representation_controller_design}, as proposed
in \cite{DiwoldKolarSchoeberl2022}. Subsequently, an additional feedback
is applied in order to stabilize reference trajectories.

In order to construct a linearizing dynamic feedback as proposed
in \cite{DiwoldKolarSchoeberl2022},
we first redefine the flat output $y_{[-R_1]}\rightarrow y$. The flat output is then of the form \eqref{eq:flat_output_disc_zta} and the corresponding parameterizing map reads as\footnote{The fact that the parameterizing map $F_{\bar{x}}$ is independent
	of $y_{[R]}$ is a consequence of the identity $\delta_{y}(F_{\bar{x}})=\tilde{f}\circ F$
	that follows from inserting the parameterizing map $F$ into \eqref{eq:implicit_euler_discretization_explicit_form}, see also \cite{DiwoldKolarSchoeberl2022}}\begin{subequations}
	\label{xxxx}
	\begin{align}
	\bar{x}&=F_{\bar{x}}(y,y_{[1]},\dots,y_{[R-1]})
	\label{state_param}\\
	u&=F_{u}(y,y_{[1]},\dots,y_{[R]})
	\label{input_param}
	\end{align}
\end{subequations}
with the multi-index $R=R_1+R_2$. Next, we extend \eqref{state_param}
by a map $z=F_{z}(y,\dots,y_{[R-1]})$ such that the combined map
$(\bar{x},z)=(F_{\bar{x}},F_{z}):=F_{\bar{x}z}$ has a regular Jacobian
matrix and is hence invertible. Since $F_{\bar{x}}$ is a submersion,
such an extension always exists. Subsequently, we define the map $\Psi(y,\dots,y_{[R]})$
given by
\begin{equation}
\label{eq:definition_psi}
\begin{aligned}\bar{x} & =F_{\bar{x}}(y,\dots,y_{[R-1]}) \\
z & =F_{z}(y,\dots,y_{[R-1]})\\
v & =y_{[R]}
\end{aligned}
\end{equation}
and its inverse $\hat{\Psi}(\bar{x},z,v)$ given by
\begin{equation}
\label{eq:definition_psi_inverse}
\begin{aligned}
(y,\dots,y_{[R-1]}) & =\hat{F}_{\bar{x}z}(\bar{x},z)\\
y_{[R]}&=v\,.
\end{aligned}
\end{equation}
As shown in \cite{DiwoldKolarSchoeberl2022}, applying the feedback
\begin{equation}
\label{eq:dynamic_feedback}
\begin{aligned}
z^{+} & =\delta_{y}(F_{z})\circ\hat{\Psi}(\bar{x},z,v)\,,\\
u&=F_{u}\circ\hat{\Psi}(\bar{x},z,v)
\end{aligned}
\end{equation}
to \eqref{eq:state_representation_controller_design} results in the
input-output behaviour $y_{[R]}=v$. Thus, in order to stabilize a
reference trajectory $y_{d}$, we can simply use an additional
feedback
\begin{equation}
v^{j}=y_{d,[r_{j}]}^{j}-\sum_{i=0}^{r_{j}-1}a_{i}^{j}(y_{[i]}^{j}\circ\hat{F}_{\bar{x}z}-y_{d,[i]}^{j})\,,\hphantom{a}j=1,\dots,m\label{eq:feedback_law}
\end{equation}
with suitable coefficients $a_{i}^{j}$, which results in a linear
tracking error dynamics. Inserting \eqref{eq:feedback_law} into \eqref{eq:dynamic_feedback}
yields a discrete-time control law of the form
\begin{align}
\begin{aligned}z^{+} & =\bar{\alpha}(\bar{x},z,y_{d},\dots,y_{d,[R]})\\
u & =\bar{\beta}(\bar{x},z,y_{d},\dots,y_{d,[R]})
\end{aligned}
\label{eq:dyn_feed}
\end{align}
which depends on the system state $\bar{x}$, the controller state
$z$, and the reference trajectory $y_{d}$ as well as its forward-shifts.
Since we assume that only the original states $x$ can be measured
or estimated, we finally have to replace $\bar{x}$ in \eqref{eq:dyn_feed}
by \eqref{eq:state_transformation} and get a control law
\begin{align}
\begin{aligned}z^{+} & ={\alpha}({x},z,y_{d},\dots,y_{d,[R]})\\
u & ={\beta}({x},z,y_{d},\dots,y_{d,[R]})\,,
\end{aligned}
\label{eq:dyn_feedx}
\end{align}
which can now be implemented on a digital processor. It should be noted again that for the design of the control law the state representation \eqref{eq:state_representation_controller_design} of the discretized system is not needed.
\section{Practical Example}
\label{sec:example}

In this section, we illustrate the results of Section \ref{sec:ImpEulerContDes} by the practical example of a VTOL aircraft, also discussed e.g. in \cite{FliessLevineMartinRouchon:1999}. The continuous-time	model of a planar VTOL aircraft in state representation reads as


	\begin{align}
	\begin{aligned}\dot{x}^{1} & =x^{4} & \dot{x}^{4} & =\tfrac{(u^{1}+u^{2})}{m}\mathrm{c}(\alpha)\mathrm{s}(x^{3})+\tfrac{(u^{1}-u^{2})}{m}\mathrm{s}(\alpha)\mathrm{c}(x^{3})\\
	\dot{x}^{2} & =x^{5} & \dot{x}^{5} & =\tfrac{(u^{1}+u^{2})}{m}\mathrm{c}(\alpha)\mathrm{c}(x^{3})+\tfrac{(u^{2}-u^{1})}{m}\mathrm{s}(\alpha)\mathrm{s}(x^{3})-g\\
	\dot{x}^{3} & =x^{6} & \dot{x}^{6} & =\tfrac{(u^{2}-u^{1})}{J}(l\mathrm{c}(\alpha)+h\mathrm{s}(\alpha))
	\end{aligned}
	\label{eq:VTOL_cont}
	\end{align}
	with $x=(q_x,q_z,q_\theta,v_{x},v_{z},\omega_\theta)$, $u=(F_{1},F_{2})$, constant
	parameters $\alpha,m,g,J$ and $h$, and the abbreviations $\mathrm{s}(\cdot)=\text{\ensuremath{\sin(\cdot)}}$
	and $\mathrm{c}(\cdot)=\text{\ensuremath{\cos(\cdot)}}$.
	It is well-known
	that
	\[
	y=(x^{1}+\varepsilon\mathrm{s}(x^{3}),x^{2}+\varepsilon\mathrm{c}(x^{3}))
	\]
	with $\varepsilon=\tfrac{J\mathrm{s}(\alpha)}{m(l\mathrm{c}(\alpha)+h\mathrm{s}(\alpha))}$
	is a flat output. With the state transformation
	\begin{align}
	\begin{aligned}\bar{x}_{4}^{1} & =x^{1}+\varepsilon\mathrm{s}(x^{3}) & \bar{x}_{3}^{1} & =x^{4}+\varepsilon x^{6}\text{\ensuremath{\mathrm{c}(x^{3})}} & \bar{x}_{2}^{1} & =x^{3}\\
	\bar{x}_{4}^{2} & =x^{2}+\varepsilon\mathrm{c}(x^{3}) & \bar{x}_{3}^{2} & =x^{5}-\varepsilon x^{6}\mathrm{s}(x^{3}) & \bar{x}_{1}^{1} & =x^{6}
	\end{aligned}
	\label{eq:vtol_state_transformation}
	\end{align}
	the system takes the form
	\begin{equation}
	\label{eq:vtol_in_strc_fx}
	\begin{aligned}
	\dot{\bar{x}}_{4}^{1} & =\bar{x}_{3}^{1}\\
	\dot{\bar{x}}_{4}^{2} & =\bar{x}_{3}^{2} \\
	\dot{\bar{x}}_{3}^{1} & =f_3^1(\bar{x},u) \\
	\dot{\bar{x}}_{3}^{2} & =f_3^2(\bar{x},u)\\
	\dot{\bar{x}}_{2}^{1} & =\bar{x}_{1}^{1}\\
	\dot{\bar{x}}_{1}^{1} & =f_1^1(\bar{x},u)\,,
	\end{aligned}
	\end{equation}
	and after an additional input transformation
	\begin{equation}
	\label{eq:input_transformation_VTOL}
	\begin{aligned}
	\bar{u}_2^1&=f_3^1(\bar{x},u) & \bar{u}_0^1&=f_1^1(\bar{x},u)
	\end{aligned}
	\end{equation}
	we get
	\begin{equation}
	\label{eq:vtol_in_strc_f}
	\begin{aligned}
	\dot{\bar{x}}_{4}^{1} & =\bar{x}_{3}^{1}\\
	\dot{\bar{x}}_{4}^{2} & =\bar{x}_{3}^{2} \\
	\dot{\bar{x}}_{3}^{1} & =\bar{u}_{2}^{1} \\
	\dot{\bar{x}}_{3}^{2} & =\tfrac{\bar{u}_{2}^{1}}{\tan(\bar{x}_{2}^{1})}-g\\
	\dot{\bar{x}}_{2}^{1} & =\bar{x}_{1}^{1}\\
	\dot{\bar{x}}_{1}^{1} & =\bar{u}_{0}^{1}\,.
	\end{aligned}
	\end{equation}
	This corresponds to the structurally flat triangular form \eqref{eq:structurally_flat_triangular_form_continuous_time_system} consisting of $p=4$ blocks with $\dim(\bar{x}_4)=2$, $\dim(\bar{x}_3)=2$, $\dim(\bar{x}_2)=1$ and $\dim(\bar{x}_1)=1$. The flat output $y$ is given by the states $\bar{x}_4$ and the states $\hat{x}_3,\dots,\hat{x}_1$ are given by $\bar{x}_3,\dots,\bar{x}_1$.
	The implicit Euler-discretization of \eqref{eq:vtol_in_strc_f} reads as
		\begin{equation}
	\begin{aligned}
	{\bar{x}}_{4}^{1,+} & =x_4^1+T_s\bar{x}_{3}^{1,+}\\
	{\bar{x}}_{4}^{2,+} & =x_4^2+T_s\bar{x}_{3}^{2,+} \\
	{\bar{x}}_{3}^{1,+} & =x_3^1+T_s\bar{u}_{2}^{1} \\
	{\bar{x}}_{3}^{2,+} & =x_3^2+T_s\tfrac{\bar{u}_{2}^{1}}{\tan(\bar{x}_{2}^{1,+})}-T_sg\\
	{\bar{x}}_{2}^{1,+} & =x_2^1+T_s\bar{x}_{1}^{1,+}\\
	{\bar{x}}_{1}^{1,+} & =x_1^1+T_s\bar{u}_{0}^{1}\,,
	\end{aligned}
	\end{equation}
	and is flat\footnote{To be precise, its explicit state representation \eqref{eq:implicit_euler_discretization_explicit_form} is flat.} with the flat output $y=\bar{x}_4$, cf. Lemma \ref{lem:preserving_flatness_of_str_form}. The corresponding parameterizing map is of the form
	\begin{equation}
	\label{eq:param_bsp_bar_bar}
	(\bar{x},\bar{u})=F(y^1_{[-3]},y^2_{[-3]},\dots,y^1_{[1]},y^2_{[1]})\,.
	\end{equation}
	 Due to Lemma \ref{lem:input_one_to_one}, the implicit Euler-discretization of the continuous-time system $\dot{\bar{x}}=f(\bar{x},u)$
	obtained from \eqref{eq:system_equation_continuous} by applying
	only the state transformation \eqref{eq:vtol_state_transformation}
	is also flat. By means of the map $u={\Phi}^{-1}_{u}(\bar{x}^{+},\bar{u})$, which is determined by \eqref{eq:input_transformation_VTOL}, the corresponding parameterization can be computed immediately from \eqref{eq:param_bsp_bar_bar} and reads as
	\begin{equation}
	\label{eq:parameterizing_map_vtol}
	(\bar{x},u)=F(y^1_{[-3]},y^2_{[-3]},\dots,y^1_{[1]},y^2_{[1]})\,.
	\end{equation}
	
	As proposed in Section \ref{subsec:Controller-Design}, based on the parameterizing map \eqref{eq:parameterizing_map_vtol} we can determine a flatness-based tracking control law systematically. For this purpose we redefine the flat output as the third backward-shift of $y$. In other words, we simply substitute $y^1_{[-3]}$ by $y^1$ and $y_{[-3]}^2$ by $y^2$, respectively. The parameterizing map reads then as
	\begin{align*}
	(\bar{x},u) & =F(y^{1},y^{2},\dots,y_{[4]}^{1},y_{[4]}^{2})\,,
	\end{align*}
	while the flat output is of the form \eqref{eq:flat_output_disc_zta}. If we choose the map $F_{z}$ as 
	\begin{equation}
	\begin{aligned}
	z^{1} & =y^{2}\\
	z^{2} & =y_{[1]}^{2}\,,
	\end{aligned}
	\end{equation}
	then the combined map $F_{\bar{x}z}$ is invertible, and a linearizing
	dynamic feedback \eqref{eq:dynamic_feedback} can be derived. The
	use of a further feedback \eqref{eq:feedback_law} results in two
	decoupled linear tracking error systems of order $(4,4)$. 
	
\begin{figure}
		\centering\includegraphics[width=1\columnwidth]{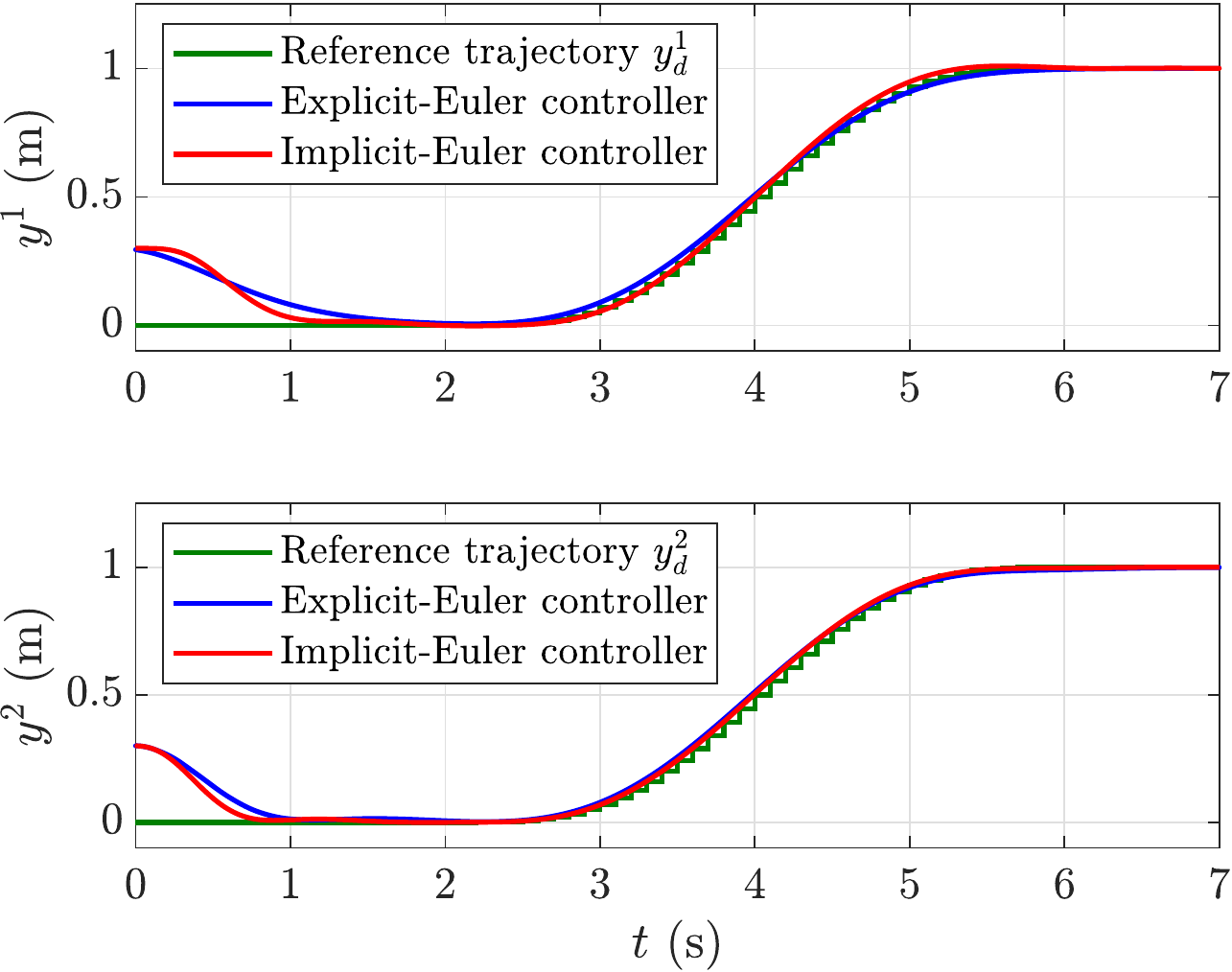}
		\caption{\label{fig:VTOL_sim_res}Tracking performance planar VTOL aircraft.}
	\end{figure}

	In a simulation
	we have compared the proposed dynamic feedback controller based on
	an implicit Euler-discretization of the VTOL with a dynamic feedback
	controller based on an explicit Euler-discretization as it is discussed
	in \cite{DiwoldKolarSchoeberl2022-1}. Both discrete-time control laws were applied to the continuous-time model \eqref{eq:VTOL_cont}, with the control inputs held piecewise constant during each sampling interval with $T_s=100\mathrm{ms}$. The step size for the numerical solver of the nonlinear differential equations \eqref{eq:VTOL_cont} was set to $T_n=0.1\mathrm{ms}$. In Fig. \ref{fig:VTOL_sim_res} 
	the simulation results are depicted, where it can be observed e.g. that the controller based on the implicit Euler-discretization
	reacts faster to the non-consistent initial conditions.

\section{Conclusion}
In this contribution we have adapted the method of \cite{DiwoldKolarSchoeberl2022-1} to use an implicit Euler-discretization instead of an explicit Euler-discretization. The main difference is that the resulting sampled-data model is flat, but not necessarily forward-flat. However, also discrete-time systems in the extended sense including backward-shifts allow for the systematic planning of trajectories and design of tracking controllers, as we have illustrated by the example of a VTOL aircraft. Subject to further research are alternative discretization methods that preserve flatness and do not necessarily require the existence of a transformation into the structurally flat triangular form \eqref{eq:structurally_flat_triangular_form_continuous_time_system} for continuous-time flat systems \eqref{eq:system_equation_continuous}.
\bibliography{Bibliography_Johannes_Mai_2022}             
                                                   






\end{document}